\documentclass[12pt,twoside]{article}
\usepackage{amssymb}
\font\teneufm=eufm10 scaled \magstep1
\font\seveneufm=eufm7 scaled \magstep1
\font\fiveeufm=eufm5  scaled \magstep1
\newfam\eufmfam
\textfont\eufmfam=\teneufm
\scriptfont\eufmfam=\seveneufm
\scriptscriptfont\eufmfam=\fiveeufm
\def\frak#1{{\fam\eufmfam\relax#1}}

\usepackage{mathrsfs}

\newfam\msbfam
\font\tenmsb=msbm10 scaled \magstep1  \textfont\msbfam=\tenmsb
\font\sevenmsb=msbm7 scaled \magstep1 \scriptfont\msbfam=\sevenmsb
\font\fivemsb=msbm5 scaled \magstep1  \scriptscriptfont\msbfam=\fivemsb
\def\Bbb{\fam\msbfam \tenmsb}

\def\RR{{\Bbb R}}
\def\CC{{\Bbb C}}

\def\QQ{{\Bbb Q}}

\def\ZZ{{\Bbb Z}}

\def\ra{\rightarrow}

 \def\HollowBoxx #1#2#3{{\dimen0=#1 \advance\dimen0 by -#2
       \dimen1=#1 \advance\dimen1 by #3
        \vrule height 0pt depth #3 width #2
       \hskip -#3
       \vrule height #1 depth #3 width #3}}
 \def\LeftContraction{\mathord{\kern1.45pt \HollowBoxx{6pt}{3.5pt}{.4pt}}\,}

 \def\HollowBox #1#2#3{{\dimen0=#1 \advance\dimen0 by -#3
       \dimen1=#1 \advance\dimen1 by #3
        \vrule height #1 depth #3 width #3
        \vrule height 0pt depth #3 width #2
        \hskip -#3}}
 \def\RightContraction{\mathord{\, \HollowBox{6pt}{3.1pt}{.4pt}} \kern1.6pt}

\def\qed{{\hfill $\Box$}}
\newtheorem{theorem}{THEOREM}[section]
\newtheorem{corollary}[theorem]{Corollary}

\newtheorem{proposition}[theorem]{Proposition}

\begin{document}

\begin{center}
{\Large \bf On Chern-Moser Normal Forms
\vspace{0.1cm}\\
of Strongly Pseudoconvex Hypersurfaces
\vspace{0.3cm}\\
with High-Dimensional Stability Group}\footnote{{\bf Mathematics Subject Classification:} 32V40, 32C05}\footnote{{\bf
Keywords and Phrases:} Chern-Moser normal forms, strongly pseudoconvex hypersurfaces, local CR-automorphisms}
\vspace{0.3cm}\\
\normalsize A. V. Isaev
\end{center}

\begin{quotation} \small \sl We explicitly describe germs of strongly pseudoconvex non-spherical real-analytic hypersurfaces $M$ at the origin in $\CC^{n+1}$ for which the group of local CR-automorphisms preserving the origin has dimension $d_0(M)$ equal to either $n^2-2n+1$ with $n\ge 2$, or $n^2-2n$ with $n\ge 3$. The description is given in terms of equations defining hypersurfaces near the origin, written in the Chern-Moser normal form. These results are motivated by the classification of locally homogeneous Levi non-degenerate hypersurfaces in $\CC^3$ with $d_0(M)=1,2$ due to A. Loboda, and complement earlier joint work by V. Ezhov and the author for the case $d_0(M)\ge n^2-2n+2$.
\end{quotation}

\thispagestyle{empty}

\pagestyle{myheadings}
\markboth{A. V. Isaev}{Hypersurfaces with High-Dimensional Stability Group}

\setcounter{section}{0}

\section{Introduction}
\setcounter{equation}{0}

Let $M$ be a strongly pseudoconvex real-analytic hypersurface in $\CC^{n+1}$ passing through the origin. In some local holomorphic coordinates $z=(z_1,\dots,z_n)$, $w=u+iv$ in a neighborhood of the origin, $M$ is given by an equation written in the Chern-Moser normal form (see \cite{CM})
\begin{equation}
v=|z|^2+\displaystyle\sum_{k,\overline{l}\ge 2}F_{k\overline{l}}(z,\overline{z},u),\label{normal}
\end{equation}
where $|z|$ is the norm of the vector $z$, and $F_{k\overline{l}}(z,\overline{z},u)$ are polynomials of
degree $k$ in $z$ and $\overline{l}$ in $\overline{z}$  whose coefficients are
analytic functions of $u$ such that the following conditions hold
\begin{equation}
\begin{array}{lll}
\hbox{tr}\,F_{2\overline{2}}&\equiv&0,\\
\hbox{tr}^2\,F_{2\overline{3}}&\equiv&0,\\
\hbox{tr}^3\,F_{3\overline{3}}&\equiv&0.
\end{array}\label{formnfcond}
\end{equation}
Here the operator $\hbox{tr}$ is defined as 
$$
\displaystyle \hbox{tr}:=\sum_{\alpha=1}^n \frac{\partial^2 }{\partial z_{\alpha}\partial \overline{z}_{\alpha}}.
$$
Everywhere below we assume that the equation of $M$ is given in the normal form.

Let $\hbox{Aut}_0(M)$ denote the {\it stability group}\, of $M$ at the origin, i.e. the group of all local CR-automorphisms of $M$ defined near the origin and preserving it. Every element $\varphi$ of $\hbox{Aut}_0(M)$ extends to a biholomorphic mapping defined in a neighborhood of the origin in $\CC^{n+1}$ and therefore can be written as
$$
\begin{array}{lll}
z&\mapsto& f_{\varphi}(z,w),\\
w&\mapsto& g_{\varphi}(z,w),
\end{array}
$$
where $f_{\varphi}$ and $g_{\varphi}$ are holomorphic. We equip $\hbox{Aut}_0(M)$ with the topology of uniform convergence of the partial derivatives of all orders of the component functions on neighborhoods of the origin in $M$. The group $\hbox{Aut}_0(M)$ with this topology is a topological group.

It is shown in \cite{CM} that every element $\varphi=(f_{\varphi},g_{\varphi})$ of $\hbox{Aut}_0(M)$ is uniquely determined by a set of parameters $(U_{\varphi}, a_{\varphi}, \lambda_{\varphi}, r_{\varphi})$, where $U_{\varphi}$ lies in the unitary group $U_n$, $a_{\varphi}\in\CC^n$, $\lambda_{\varphi}>0$, $r_{\varphi}\in\RR$. These parameters are found from the following relations
$$
\begin{array}{ll}
\displaystyle\frac{\partial f_{\varphi}}{\partial z}(0)=\lambda_{\varphi}U_{\varphi}, &
\displaystyle\frac{\partial f_{\varphi}}{\partial w}(0)=\lambda_{\varphi}U_{\varphi}a_{\varphi},\\
\vspace{0.2cm} &\\
\displaystyle\frac{\partial g_{\varphi}}{\partial w}(0)=\lambda_{\varphi}^2, &
\hbox{Re}\,\displaystyle\frac{\partial^2 g_{\varphi}}{\partial^2 w}(0)=2\lambda_{\varphi}^2r_{\varphi}.
\end{array}
$$
For results on the dependence of local CR-mappings on their jets in more general settings see e.g. \cite{BER1}, \cite{BER2}, \cite{Eb}, \cite{Z}.   

We assume that $M$ is non-spherical at the origin, i.e. $M$ in a neighborhood of the origin is not CR-equivalent to an open subset of the sphere $S^{2n+1}\subset\CC^{n+1}$. In this case for every element $\varphi=(f_{\varphi},g_{\varphi})$ of  $\hbox{Aut}_0(M)$ we have $\lambda_{\varphi}=1$ and the parameters $a_{\varphi}, r_{\varphi}$ are uniquely determined by the matrix $U_{\varphi}$; moreover, the mapping 
$$
\Phi:\, \hbox{Aut}_0(M)\ra GL_n(\CC),\qquad \Phi:\,\varphi\mapsto U_{\varphi}
$$
is a topological group isomorphism between $\hbox{Aut}_0(M)$ and\linebreak $G_0(M):=\Phi(\hbox{Aut}_0(M))$, with $G_0(M)$ being a real algebraic subgroup of $GL_n(\CC)$ (see \cite{CM}, \cite{B}, \cite{L1}, \cite{BV}, \cite{VK}). We pull back the Lie group structure from $G_0(M)$ to $\hbox{Aut}_0(M)$ by means of $\Phi$ and denote by $d_0(M)$ the dimension of $\hbox{Aut}_0(M)$. Clearly, $d_0(M)\le n^2$. 

We are interested in characterizing hypersurfaces for which $d_0(M)$ is large (certainly positive). We show that in some normal coordinates the equations of such hypersurfaces take a very special form. As will be explained below, results of this kind potentially can be applied to the  classification problem for locally CR-homogeneous strongly pseudoconvex hypersurfaces. For $n=1$ this problem was solved by E. Cartan in \cite{C}. For $n=2$ with $d_0(M)>0$ an explicit classification was obtained in \cite{L3}, \cite{L5}. For $n\ge 3$ there is no such classification even for hypersurfaces with high-dimensional stability group. Note, however, that {\it globally}\, homogeneous hypersurfaces have been extensively studied (see e.g. \cite{AHR} and references therein). We also mention that locally homogeneous hypersurfaces in $\CC^3$ with non-degenerate indefinite Levi form and 2-dimensional stability group were classified in \cite{L4} and that recently in \cite{FK} Fels and Kaup have determined all locally homogeneous 5-dimensional CR-manifolds with certain degenerate Levi forms.     

For a non-spherical hypersurface $M$ the group $\hbox{Aut}_0(M)$ is known to be linearizable, i.e. in some normal coordinates every $\varphi\in\hbox{Aut}_0(M)$ can be written in the form
$$
\begin{array}{lll}
z&\mapsto& U_{\varphi}z,\\
w&\mapsto& w,
\end{array}
$$
(see \cite{KL}). If all elements of $\hbox{Aut}_0(M)$ in some coordinates have the above form, we say that $\hbox{Aut}_0(M)$ is linear in these coordinates. Thus, in order to describe hypersurfaces $M$ with a particular value of $d_0(M)$, one needs to: (a) write $M$ in normal coordinates in which $\hbox{Aut}_0(M)$ is linear, (b) determine all closed subgroups $H$ of $U_n$ of dimension $d_0(M)$, and (c) find all $H$-invariant real-analytic functions of $z$, $\overline{z}$ and $u$, homogeneous of fixed degrees in each of $z$ and $\overline{z}$. Then every $F_{k\overline{l}}(z,\overline{z},u)$ in (\ref{normal}) is a function of the kind found in (c), and one obtains the general form of $M$.

In \cite{EI} we considered the case $d_0(M)\ge n^2-2n+2$ for $n\ge 2$. It turned out that if $d_0(M)\ge n^2-2n+3$, then $d_0(M)=n^2$, that is, $G_0(M)=U_n$. Clearly, every $U_n$-invariant real-analytic function is a function of $|z|^2$ and $u$, and thus the equation of $M$ in any normal coordinates in which $\hbox{Aut}_0(M)$ is linear has the form   
\begin{equation}
v=|z|^2+\displaystyle\sum_{p=4}^{\infty}C_p(u)|z|^{2p}\label{form2}
\end{equation}
where $C_p(u)$ are real-valued analytic functions of $u$, and for some $p$ we have $C_p(u)\not\equiv 0$. Here the condition $p\ge 4$ comes from identities (\ref{formnfcond}).

Further, for $d_0(M)=n^2-2n+2$ we showed that the equation of $M$ in some normal coordinates in which $\hbox{Aut}_0(M)$ is linear has the form
\begin{equation}
v=|z|^2+\sum_{p+q\ge 2}C_{pq}(u)|z_1|^{2p}|z|^{2q},\label{form1}
\end{equation}
where $C_{pq}(u)$ are real-valued analytic functions of $u$, $C_{pq}(u)\not\equiv 0$ for some $p,q$ with $p>0$, and $C_{pq}$ for $p+q=2,3$ satisfy certain conditions arising from identities (\ref{formnfcond}).\footnote{In \cite{EI} we erroneously stated that identities (\ref{formnfcond}) imply that $C_{pq}=0$ for $p+q=2,3$. This is in general not the case (see the erratum to \cite{EI}).} Equation (\ref{form1}) is the most general form of a hypersurface with $d_0(M)=n^2-2n+2$ and cannot be simplified any further without additional assumptions on $M$. This equation is a consequence of our description of closed connected subgroups of $U_n$ of dimension $n^2-2n+2$ obtained earlier in \cite{IK}.

In \cite{L3} A. Loboda classified strongly pseudoconvex locally\linebreak CR-homogeneous hypersurfaces in $\CC^3$ with $d_0(M)=2$ (here $n=2$) by means of normal form techniques (see also \cite{L4}). Using the homogeneity of $M$ and the condition $d_0(M)=2$ he showed that the equation of $M$ must significantly simplify, which eventually yielded the classification. His arguments avoid using the explicit form of closed connected 2-dimensional subgroups of $U_2$ (every such subgroup is conjugate to $U_1\times U_1$) and, as a result, special normal form (\ref{form1}). It seems that (\ref{form1}) can be utilized to simplify the proof of the main result of \cite{L3}. Further, equation (\ref{form1}) may be a useful tool for describing locally CR-homogeneous strongly pseudoconvex hypersurfaces with $d_0(M)=n^2-2n+2$ for arbitrary $n\ge 2$. Overall, the introduction of algebraic arguments into the analysis of normal forms seems to be a fruitful idea. 

Observe for comparison that every locally CR-homogeneous strongly pseudoconvex hypersurface with $d_0(M)\ge n^2-2n+3$ and $n\ge 2$ is spherical, since by (\ref{form2}) the origin is an umbilic point of $M$. This is in contrast with hypersurfaces whose Levi form is non-degenerate and indefinite (see \cite{EI} for a description of such hypersurfaces with $d_0(M)\ge n^2-2n+3$ and \cite{L2} for the homogeneous case with $n=2$).

In this paper we consider the cases $d_0(M)=n^2-2n+1$ with $n\ge 2$, and $d_0(M)=n^2-2n$ with $n\ge 3$. Our result is the following theorem.

\begin{theorem}\label{main}\sl Let $M$ be a strongly pseudoconvex real-analytic\linebreak non-spherical hypersurface in $\CC^{n+1}$ passing through the origin.
\vspace{0.3cm}

\noindent (A) If $d_0(M)=n^2-2n+1$ and $n\ge 2$, then in some normal coordinates near the origin in which $\hbox{Aut}_0(M)$ is linear the equation of $M$ takes one of the following three forms: 
\vspace{0.3cm}

\begin{equation}
v=|z|^2+\sum_{\scriptsize\begin{array}{l}p+q\ge 2,\, r+s\ge 2,\\
(p-r)k_1+(q-s)k_2=0
\end{array}}C_{pqrs}(u)z_1^pz_2^q\overline{z_1}^r\overline{z_2}^s,\label{a3}
\end{equation}
where $k_1,k_2$ are non-zero integers with $(k_1,k_2)=1$ and $k_2>0$, $C_{pqrs}(u)$ are real-analytic functions of $u$, and $C_{pqrs}(u)\not\equiv 0$ for some $p,q,r,s$ with either $p\ne r$ or $q\ne s$ (here $n=2$);
\vspace{0.3cm}

\begin{equation}
 v=|z|^2+\sum_{2p+q\ge 2}C_{pq}(u)|z_1^2+z_2^2+z_3^2|^{2p}|z|^{2q},\label{a1}
\end{equation}
where $C_{pq}(u)$ are real-valued analytic functions of $u$, and $C_{pq}(u)\not\equiv 0$ for some $p,q$ with $p>0$ (here $n=3$);
\vspace{0.3cm}

\begin{equation}
v=|z|^2+\sum_{p+r,\, q+r\ge 2}C_{pqr}(u)z_1^p\overline{z_1}^q|z|^{2r},\label{a2}
\end{equation}
where $C_{pqr}(u)$ are real-analytic functions of $u$, and $C_{pqr}(u)\not\equiv 0$ for some $p,q,r$ with $p\ne q$.
\vspace{0.3cm}

\noindent (B) If $d_0(M)=n^2-2n$ and $n\ge 3$, then in some normal coordinates near the origin in which $\hbox{Aut}_0(M)$ is linear the equation of $M$ takes one of the following three forms: 
\vspace{0.3cm}

\begin{equation}
 v=|z|^2+\sum_{2p+r\ge 2,\,2q+r\ge 2}C_{pqr}(u)(z_1^2+z_2^2+z_3^2)^p(\overline{z_1}^2+\overline{z_2}^2+\overline{z_3}^2)^q|z|^{2r},\label{b3}
\end{equation}
where $C_{pqr}(u)$ are real-analytic functions of $u$, and $C_{pqr}(u)\not\equiv 0$ for some $p,q,r$ with $p\ne q$ (here $n=3$);
\vspace{0.3cm}

\begin{equation}
 v=|z|^2+\sum_{p+q+r\ge 2}C_{pqr}(u)|z_1|^{2p}|z_2|^{2q}|z_3|^{2r},\label{b1}
\end{equation}
where $C_{pqr}(u)$ are real-valued analytic functions of $u$, and $C_{pqr}(u)\not\equiv 0$ for some $p,q,r$ (here $n=3$);
\vspace{0.3cm}

\begin{equation}
 v=|z|^2+\sum_{p+q\ge 2}C_{pq}(u)|z'|^{2p}|z''|^{2q},\label{b2}
\end{equation}
where $z':=(z_1,z_2)$, $z'':=(z_3,z_4)$, $C_{pq}(u)$ are real-valued analytic functions of $u$, and $C_{pq}(u)\not\equiv 0$ for some $p,q$ (here $n=4$).
\end{theorem}

\begin{corollary}\label{cor}\sl Let $M$ be a strongly pseudoconvex real-analytic\linebreak non-spherical hypersurface in $\CC^{n+1}$ passing through the origin. Assume that $n\ge 5$ and $d_0(M)\ge n^2-2n$. Then $d_0(M)\ge n^2-2n+1$. Furthermore, in some normal coordinates near the origin in which $\hbox{Aut}_0(M)$ is linear the equation of $M$ has the form
$$
v=|z|^2+\sum_{p+r,\, q+r\ge 2}C_{pqr}(u)z_1^p\overline{z_1}^q|z|^{2r},
$$
where $C_{pqr}(u)$ are real-analytic functions of $u$, and $C_{pqr}(u)\not\equiv 0$ for some $p,q,r$. 
\end{corollary} 

Locally CR-homogeneous hypersurfaces in $\CC^3$ with $d_0(M)=1$ (here $n=2$) were classified in \cite{L5} and we believe that Part (A) of Theorem \ref{main} can be used to simplify Loboda's arguments.

\section{Proof of Theorem \ref{main}} 
\setcounter{equation}{0}

The main ingredient of the proof of Theorem \ref{main} is the following proposition.  

\begin{proposition}\label{un1} \sl Let $H$ be a connected closed subgroup of $U_n$ with $n\ge 2$.
\vspace{0.3cm} 

\noindent If $\hbox{dim}\,H=n^2-2n+1$, then $H$ is conjugate in $U_n$ to one of the following subgroups:
\vspace{0.3cm}

\noindent (i) $e^{i\RR}SO_3(\RR)$ (here $n=3$); 
\vspace{0.3cm}

\noindent (ii) $U_1\times SU_{n-1}$ realized as the subgroup of all matrices
$$
\left(\begin{array}{cc}
e^{i\theta} & 0\\
0& A
\end{array}\right),
$$
where $\theta\in\RR$ and $A \in SU_{n-1}$, for $n\ge 3$; 
\vspace{0.3cm}

\noindent (iii) the subgroup $H_{k_1,k_2}^n$ of all matrices
\begin{equation}
\left(\begin{array}{cc}
a & 0\\
0 & A
\end{array}\right),\label{mat}
\end{equation}
where $k_1,k_2$ are fixed integers such that $(k_1,k_2)=1$, $k_2>0$, and
$A\in U_{n-1}$, $a\in
(\det A)^{\frac{k_1}{k_2}}:=\exp(k_1/k_2\, \hbox{Ln}\,(\det A))$.\footnote{For $k_1\ne 0$ the group $H_{k_1,k_2}^n$ is a $k_2$-sheeted cover of $U_{n-1}$.}
\vspace{0.3cm}

\noindent If $\hbox{dim}\,H=n^2-2n$, then $H$ is conjugate in $U_n$ to one of the following subgroups:
\vspace{0.3cm}

\noindent (iv) $SO_3(\RR)$ (here $n=3$); 
\vspace{0.3cm}

\noindent (v) $U_1\times U_1\times U_1$ realized as diagonal matrices in $U_3$ (here $n=3$); 
\vspace{0.3cm}

\noindent (vi) $U_2\times U_2$ realized as block-diagonal matrices in $U_4$ (here $n=4$); 
\vspace{0.3cm}

\noindent (vii) $SU_{n-1}$ realized as the subgroup of all matrices
$$
\left(\begin{array}{cc}
1 & 0\\
0& A
\end{array}\right),\quad A \in SU_{n-1}. 
$$
\end{proposition}

\noindent {\bf Proof:} Suppose first that $\hbox{dim}\,H=n^2-2n+1$. Since $H$ is compact, it is completely
reducible, i.e. $\CC^n$ splits into the sum of $H$-invariant
pairwise orthogonal complex subspaces, $\CC^n=V_1\oplus\dots\oplus V_m$,
such that the restriction $H_j$ of $H$ to each  $V_j$ is irreducible. Let
$n_j:=\hbox{dim}_{\CC}V_j$ (hence $n_1+\dots+n_m=n$) and let
$U_{n_j}$ be the group of unitary
transformations of $V_j$. Clearly, $H_j\subset U_{n_j}$, and therefore
$\hbox{dim}\,H\le n_1^2+\dots+n_m^2$.
On the other hand  $\hbox{dim}\,H=n^2-2n+1$, which shows that
$m\le 2$.

Let $m=2$. Then there exists a unitary change of coordinates
in $\CC^n$ such all elements of $H$ take the form (\ref{mat}), where $A \in U_{n-1}$ and $a \in U_1$. If $\hbox{dim}\,H_1=0$, then $H_1=\{1\}$, and
therefore $H_2=U_{n-1}$. In this case we obtain the group $H_{0,1}^n$. Suppose next that $\hbox{dim}\,H_1=1$, i.e. $H_1=U_1$. Then
$n^2-2n\le\hbox{dim}\,H_2\le n^2-2n+1$. If $\hbox{dim}\,H_2=n^2-2n$, then $H_2=SU_{n-1}$, and hence $H$ is conjugate to
$U_1\times SU_{n-1}$ for $n\ge 3$ and to $H_{0,1}^2$ for $n=2$. Now let  $\hbox{dim}\,H_2=n^2-2n+1$, i.e. $H_2=U_{n-1}$. Consider the Lie algebra ${\frak h}$ of $H$. Up to conjugation, it consists of matrices of the form
\begin{equation}
\left(\begin{array}{cc}
l({\frak A}) & 0\\
0& {\frak A}
\end{array}\right),\label{mat2}
\end{equation}
where ${\frak A}\in{\frak u}_{n-1}$ and $l({\frak A})\not\equiv 0$
is a linear function  of the matrix elements of ${\frak A}$ ranging  in $i\RR$. Clearly, $l({\frak A})$ must vanish on the derived algebra of ${\frak u}_{n-1}$, which is ${\frak{su}}_{n-1}$. Hence matrices (\ref{mat2}) form  a Lie algebra if and only if $l({\frak A})=c\cdot\hbox{trace}\,{\frak A}$, where $c\in\RR\setminus\{0\}$. Such an algebra can be the Lie algebra of a closed subgroup of $U_{n-1}\times U_1$  only if $c\in\QQ\setminus\{0\}$. Hence $H$ is conjugate to $H_{k_1,k_2}^n$ for some $k_1,k_2\in\ZZ$, where one can always assume that $k_2>0$ and $(k_1,k_2)=1$.

Now let $m=1$. We  shall proceed as in the proof of Lemma 1.4 in \cite{I}. Let ${\frak h}^{\CC}:={\frak h}+i{\frak h}\subset{\frak{gl}}_n$ be the complexification of ${\frak h}$, where ${\frak{gl}}_n:={\frak{gl}}_n(\CC)$. The algebra ${\frak h}^{\CC}$ acts irreducibly on $\CC^n$ and by a theorem of E. Cartan, ${\frak h}^{\CC}$ is either semisimple or  the direct sum of the center ${\frak c}$ of ${\frak{gl}}_n$ and a semisimple ideal ${\frak t}$. Clearly, the action of the ideal ${\frak t}$ on $\CC^n$ is irreducible.

Assume first that ${\frak h}^{\CC}$ is semisimple, and let ${\frak
h}^{\CC}={\frak h}_1\oplus\dots\oplus{\frak h}_k$ be its decomposition
into the direct sum of simple ideals. Then the natural irreducible $n$-dimensional representation of
${\frak h}^{\CC}$ (given by the embedding of ${\frak h}^{\CC}$ in ${\frak{gl}}_n$) is the tensor product of some irreducible faithful representations of the ${\frak h}_j$. Let $n_j$ be the dimension of the corresponding representation  of ${\frak h}_j$, $j=1,\dots,k$. Then $n_j\ge 2$, $\hbox{dim}_{\CC}\,{\frak h}_j\le n_j^2-1$, and $n=n_1\cdot...\cdot n_k$. 

It is straightforward to show that if $n=n_1\cdot...\cdot n_k$ with $k\ge 2$ and $n_j\ge 2$ for $j=1,\dots,k$, then $\sum_{j=1}^k n_j^2\le n^2-2n$. Since $\hbox{dim}_{\CC}\,{\frak h}^{\CC}=n^2-2n+1$, it then follows that $k=1$, i.e. ${\frak h}^{\CC}$ is
simple. The minimal dimensions of irreducible faithful
representations $V$ of complex simple Lie algebras ${\frak s}$ are well-known and shown in the following table (see e.g. \cite{OV}).

\begin{center}
\begin{tabular}{|l|c|c|}
\hline
\multicolumn{1}{|c|}{${\frak s}$}&
\multicolumn{1}{c|}{ $\hbox{dim}\,V$}&
\multicolumn{1}{c|}{$\hbox{dim}\,{\frak s}$}
\\ \hline
${\frak{sl}}_k$\,\,$k\ge 2$ & $k$ & $k^2-1$
\\ \hline
${\frak o}_k$\,\, $k\ge 7$&  $k$ & $k(k-1)/2$
\\ \hline
${\frak{sp}}_{2k}$\,\,$k\ge 2$ & $2k$ & $2k^2+k$
\\ \hline
${\frak e}_6$ & 27 & 78
\\ \hline
${\frak e}_7$ & 56 & 133
\\ \hline
${\frak e}_8$ & 248 & 248
\\ \hline
${\frak f}_4$ & 26 & 52
\\ \hline
${\frak g}_2$ & 7 & 14
\\ \hline
\end{tabular}
\end{center}
It is straightforward to see that none of these dimensions is compatible with the condition $\hbox{dim}_{\CC}\,{\frak h}^{\CC}=n^2-2n+1$. Therefore, ${\frak h}^{\CC}={\frak c}\oplus{\frak t}$, where $\hbox{dim}\,{\frak t}=n^2-2n$. Then, if $n=2$, we obtain that $H$ coincides with the center of $U_2$ which is impossible since its action on $\CC^2$ is then not irreducible. Assuming that $n\ge 3$ and repeating the above argument for ${\frak t}$ in place of ${\frak h}^{\CC}$, we see that ${\frak t}$ can only be isomorphic to ${\frak{sl}}_{n-1}$. But ${\frak{sl}}_{n-1}$ does not have an irreducible $n$-dimensional representation unless $n=3$. 

Thus, $n=3$ and ${\frak h}^{\CC}\simeq\CC\oplus{\frak{sl}}_2\simeq\CC\oplus{\frak{so}}_3$. Further, we observe that every irreducible $3$-dimensional representation of ${\frak{so}}_3$ is equivalent to its defining representation. This implies that $H$ is conjugate in $GL_3(\CC)$ to $e^{i\RR}SO_3(\RR)$. Since $H\subset U_3$ it is straightforward to show that the conjugating element can be chosen to belong to $U_3$. This completes the proof of the proposition in the case $\hbox{dim}\,H=n^2-2n+1$.

For $\hbox{dim}\,H=n^2-2n$ we argue analogously and see that either $m\le 2$, or, for $n=3$ we have $m=3$. In the latter case $H$ is conjugate in $U_3$ to $U_1\times U_1\times U_1$.

Let $m=2$. Then either $n=4$ and $H$ is conjugate in $U_4$ to $U_2\times U_2$, or there exists a unitary change of coordinates
in $\CC^n$ such all elements of $H$ take the form (\ref{mat}), where $A \in U_{n-1}$ and $a \in U_1$. If $\hbox{dim}\,H_1=0$, then $H_1=\{1\}$, and therefore $H_2=SU_{n-1}$. Assume now that $\hbox{dim}\,H_1=1$, i.e. $H_1=U_1$. Then $n\ge 3$ and $n^2-2n-1\le\hbox{dim}\,H_2\le n^2-2n$. Lemma 1.4 of \cite{I} shows that the possibility $\hbox{dim}\,H_2=n^2-2n-1$ cannot in fact occur, and thus we have $\hbox{dim}\,H_2=n^2-2n$. Then $H_2=SU_{n-1}$, and hence $H$ is conjugate to a codimension 1 subgroup of the group of all matrices of the form (\ref{mat}) with $A\in SU_{n-1}$. Consider the Lie algebra ${\frak h}$ of $H$. Up to conjugation, it consists of matrices of the form (\ref{mat2}), where ${\frak A}\in{\frak{su}}_{n-1}$ and $l({\frak A})\not\equiv 0$ is a linear function  of the matrix elements of ${\frak A}$ ranging  in $i\RR$. Clearly, $l({\frak A})$ must vanish on the derived algebra of ${\frak{su}}_{n-1}$, which is ${\frak{su}}_{n-1}$ itself. This contradiction shows that the possibility $\hbox{dim}\,H_1=1$ does not in fact realize. 

In the case $m=1$ we argue as in the case $\hbox{dim}\,H=n^2-2n+1$. If ${\frak h}^{\CC}$ is semisimple, it follows as before that ${\frak h}^{\CC}$ is in fact simple. A glance at the table of minimal dimensions of irreducible faithful representations of complex simple Lie algebras now yields that $n=3$ and ${\frak h}^{\CC}\simeq{\frak{sl}}_2\simeq{\frak{so}}_3$, and hence $H$ is conjugate in $U_3$ to $SO_3(\RR)$. If, finally, $n\ge 3$ and ${\frak h}^{\CC}={\frak c}\oplus{\frak t}$, where $\hbox{dim}\,{\frak t}=n^2-2n-1$, we see that ${\frak t}$ must be simple and obtain a contradiction with the above table.

The proof of the proposition is complete.\qed
\smallskip\\

To finalize the proof of Theorem \ref{main} we now need to determine polynomials in $z$, $\overline{z}$ with coefficients depending on $u$, invariant under each of the groups listed in (i)-(vii) of Proposition \ref{un1}. This is not hard to do. Indeed, every $SO_3(\RR)$-invariant polynomial is a function of $z_1^2+z_2^2+z_3^2$, $\overline{z_1}^2+\overline{z_2}^2+\overline{z_3}^2$ and $|z|^2$. If, in addition, such a polynomial is $e^{i\RR}$-invariant, it depends only on $|z_1^2+z_2^2+z_3^2|^2$ and $|z|^2$. These observations lead to forms (\ref{a1}) and (\ref{b3}). Next, $U_1\times SU_{n-1}$-invariant polynomials for $n\ge 3$ are in fact $U_1\times U_{n-1}$-invariant and therefore lead to hypersurfaces with $d_0(M)\ge n^2-2n+2$. Further, every $H_{0,1}^n$-invariant polynomial is a function of $z_1$, $\overline{z_1}$ and $|z|^2$, which leads to form (\ref{a2}). Every $H_{k_1,k_2}^n$-invariant polynomial for $k_1\ne 0$ and $n\ge 3$ is in fact $U_1\times U_{n-1}$-invariant; such polynomials lead to hypersurfaces with $d_0(M)\ge n^2-2n+2$. Observe also that invariance under the group $H_{k_1,k_2}^2$ with $k_1\ne 0$ (here $n=2$) leads to form (\ref{a3}). Next, $U_1\times U_1\times U_1$-invariant polynomials are functions of $|z_1|^2$, $|z_2|^2$, $|z_3|^3$ and lead to form (\ref{b1}). Similarly, $U_2\times U_2$-invariant polynomials are functions of $|z'|^2$, $|z''|^2$, where $z':=(z_1,z_2)$, $z'':=(z_3,z_4)$, and lead to form (\ref{b2}). Finally, $SU_{n-1}$-invariant polynomials for $n\ge 3$ are in fact $U_{n-1}$-invariant and hence lead to hypersurfaces with $d_0(M)\ge n^2-2n+1$.

The proof of Theorem \ref{main} is complete.\qed

{\obeylines
Department of Mathematics
The Australian National University
Canberra, ACT 0200
AUSTRALIA
E-mail: alexander.isaev@maths.anu.edu.au
}

\end{document}